\documentclass[12pt,fleqn]{article}
\usepackage{amsmath}
\usepackage{amsfonts}
\usepackage{amsthm}

\newcommand{\Rset}{\mathbb{R}}
\newcommand{\Cset}{\mathbb{C}}

\begin{document}

\title{Random integral representations for free-infinitely divisible and tempered
stable distributions\footnote{Research funded in part by a grant
Nr 1P03A04629, 2005-2008}}

\author{Zbigniew J. Jurek}

\date{In: \textbf{Stat. Probab. Letters 77(2007), pp. 417-425.}}

\maketitle

\newtheorem{prop}{PROPOSITION}
\newtheorem{cor}{COROLLARY}
\newtheorem{thm}{Theorem}
\theoremstyle{remark}
\newtheorem{rem}{REMARK}

\begin{quote}
\noindent {\footnotesize \textbf{ABSTRACT.}  There  are given
sufficient conditions under which mixtures of dilations of L\'evy
spectral measures, on a Hilbert space, are L\'evy measures again.
We introduce some random integrals with respect to infinite
dimensional L\'evy
 processes, which in turn give some integral mappings. New classes
(convolution semigroups) are introduced. One of them gives an
unexpected relation between the free (Voiculescu) and the
classical L\'evy-Khintchine formulae while the second one
coincides with tempered stable measures (Mantegna nad Stanley)
arisen in statistical physics.}

\medskip
MSC\,2000 \emph{subject classifications.} Primary 60E07, 60E10 ;
secondary 60B11, 60G51.

\end{quote}

\newpage
In probability theory and mathematical statistics the Fourier and
Laplace transforms are the main tools used to prove weak limit
theorems or to identify probability distributions. These are
purely \emph{analytic methods} used in other branches of
mathematics as well. Jurek and Vervaat (1983) had introduced
\emph{the method of random integral representations} that allows
to describe distributions as laws of some random integrals with
respect to L\'evy processes. In fact, this method was introduced
for a study of the class $L$ of selfdecomposable probability
measures on Banach spaces. Namely, it was proved that for a
probability measure $\mu$, on a Banach space $E$, we have
\begin{equation}
\mu\in L \ \mbox{iff} \ \
\mu=\mathcal{L}(\int_0^{\infty}\,e^{-t}\,dY(t)) \ \ \mbox{for a
unique L\'evy process $Y$.}
\end{equation}
More precisely, to insure the existence of the above improper
random integrals one needs that
$\mathbb{E}[\log(1+||Y(1)||)]<\infty$, and then the induced random
integral mapping \emph{$\mathcal{I}$}
\[
ID_{\log}(E)\ni
\mathcal{L}(Y(1))\to\mathcal{L}(\int_0^{\infty}\,e^{-t}\,dY(t))\in
L(E)
\]
is an isomorphism between the corresponding convolution semigroups
of probability measures: infinitely divisible with finite
logarithmic moments $ID_{\log}$ and selfdecomposable ones $L$. To
the (unique) L\'evy process $Y$ in (1) one refers as the
\emph{background driving L\'evy process} for the selfdecomposable
measure $\mu$, in short its BDLP. Besides "randomness" the
\emph{stochastic method} given by the integral representation (1),
gives easily characterization of measures in terms of their
Fourier transforms. Moreover, it allows to incorporate space and
time changes in order to retrieve the $BDLP$, i.e.,
\[
\mathcal{L}(\int_0^{c}\,e^{-t}\,dY(t/c))\Rightarrow
\mathcal{L}(Y(1)), \ \ \mbox{as} \ \ c\to 0.
\]

The aim of this note is to extend the random integral
representation approach to more general schemes then those in
Jurek (1982, 1983, 1985, 1988) and Jurek and Vervat (1983). As a
consequence we will find such representations for two classes
(semigropus) of measures. First, a class $\mathcal{E}$ that
coincides with a class of free infinitely divisible measures in
Voiculescu sense (cf. Barndorff-Nielsen and Thorbjornsen (2002))
and the second, $\mathcal{TS}$ coincides with the tempered stable
distributions (these are a generalization of the titled stable
measures) that arose among others in statistical physics; (cf.
Mantegna and Stanley (1994), Novikov (1995), Kaponen (1995) and
Rosi\'nski (2002), (2004)). One may expect that the random
integral representations will help in some simulation problems.

Our results are mainly given in a generality of measures on a
Hilbert space, with some digression to a case of Banach spaces.
However, methods of proofs are dimensionless so they can be read
for Euclidean spaces as well.

\medskip
 \textbf{1. Introduction and terminology.} Let $E$
denotes a real separable Banach space, $E^{\prime}$ its conjugate
space, $< \cdot,\cdot >$ the usual pairing between $E$ and
$E^{\prime}$, (this is just a scalar product in a case  when $E$
is a Hilbert or an Euclidean space), and $||.||$ the norm on E.
The $\sigma$-field of all Borel subsets of $E$ is denoted by
$\mathcal{B}$, while $\mathcal{B}_{0}$ denotes Borel subsets of $E
\setminus \{0\}$. By $\mathcal{P}(E)$ or simply by $\mathcal{P}$,
we denote the (topological) semigroup of all Borel probability
measures on $E$, with convolution ``$\ast$" and the topology of
weak convergence ``$\Rightarrow$."

Recall that a measure $\mu\in \mathcal{P}$ is called
\emph{infinitely divisible} if for each natural  $n\ge 2$ there
exists $\nu_n \in \mathcal{P}$ such that $\nu_n^{\star n} = \mu$.
The class $ID(E)$, of all infinitely divisible probability
measures on $E$, is a closed topological  subsemigroup of
$\mathcal{P}$. Each ID distribution $\mu$ is uniquely determined
by a triple: a shift vector $a\in E$, a Gaussian covariance
operator $R$, and a L\'evy spectral measure $M$; we will write
$\mu = [a,R,M]$. These are the parameters in the L\'evy-Khintchine
representation of the characteristic function $\hat{\mu}$, namely
\begin{multline}
\mu \in ID \quad \mbox{iff} \quad \hat{\mu}(y)= \exp[\,\Phi(y)], \
\ \
\mbox{where}  \\
\Phi(y):= i<y,a>-1/2<Ry,y>+ \\
\int_{E \setminus\{0\}} [e^{i<y,x>}-1-i<y,x>1_{||x||\leq
1}(x)]M(dx),\ \ y \in E^{\prime};
\end{multline}
$\Phi$ is called the \textit{L\'evy exponent} of $\hat{\mu}$ (cf.
Araujo and Gin\'e (1980), Section 3.6). For $\mu\in ID(E)$ one can
define arbitrary positive convolution powers, namely
\begin{equation}
\mbox{if} \ \mu=[a, R, M] \in ID(E) \ \mbox{then} \ \mu^{\star
\,c}=[c\cdot a, c\cdot R, c\cdot M], \ \mbox{for any} \ c\ge 0,
\end{equation}
and
\begin{equation}
 T_a\mu =[a_c,c^2R,T_cM], \mbox{where} \ a_c = c\,a + c\int_E
 x[1_B(cx)-1_B(x)]M(dx)
\end{equation}
where the mapping $T_c$ is given by $T_cx=cx, x\in E$ and
$T_cM(A)=M(c^{-1}A)$ for all $A \in \mathcal{B}_{0}$.

In a case when $\rho$ is the probability distribution of an
$E$-valued random element $\xi$ then $T_a\rho$ is the probability
distribution of the random element $a \cdot \xi$, for any real
number $a$.

Let $\mathcal{M}(E)$ denotes the totality of all L\'evy spectral
measures on a space $E$. It is a positive cone but it is also
closed under dilations $T_a$, i.e.,
\begin{equation}
M\in \mathcal{M}(E) \ \ \mbox{iff} \ \ T_a\,M \in \mathcal{M}(E),
\ \mbox{for all} \ \ a\in \Rset.
\end{equation}
Hence we conclude that
\begin{equation}
M\in \mathcal{M}(E) \ \mbox{iff}\ \sum_{i=1}^k\,c_i\cdot
T_{a_i}\,M \in\mathcal{M}(E) \ \mbox{for all}\ \ k \ge 1, c_i >0,
a_i \in \Rset.
\end{equation}

\medskip
\medskip
\textbf{2. The $\lambda$ -mixtures of L\'evy spectral measures.}
 Our main objective in this section is to provide a method of
constructing L\'evy spectral measures using the mixtures of
$T_tM$. Namely, for a non-negative Borel measure $\lambda$ on
$\Rset^+=(0, \infty)$ and a Borel measure $M$ on $E\setminus\{0\}$
we define
\begin{equation}
M^{(\lambda)}(A):=\int_0^{\infty}\,(T_tM)(A)\lambda(dt)=
\int_0^{\infty}\int_{E\setminus\{0\}}\,1_A(tx)M(dx)\lambda(dt),\ \
A\in\mathcal{B}_{0},
\end{equation}
and call it the \emph{$\lambda$-mixture} of $T_tM, t>0$. Questions
on mixtures of L\'evy spectral measures on Banach spaces were
investigated in Jurek (1990). In particular, it was proved that
\begin{multline}
\mbox{\emph{if}} \ \ M^{(\lambda)}\in \mathcal{M}(E) \ \
\mbox{\emph{then}} \ \ M\in \mathcal{M}(E), \ \  \int_0^{\infty}
\min(1, t^2)\,\lambda(dt)<\infty \ \ \mbox{\emph{and}}\\
\int_0^{\infty}M(\{x:
||x||>t^{-1}\})\lambda(dt)=\int_{E\setminus\{0\}}\lambda(\{s:s>||x||^{-1}\})M(dx)<\infty;
\end{multline}
cf. Jurek (1990), Proposition 2. The converse implication to (8)
is not completely settled. Note that for measures $\lambda$ with
finite support and arbitrary L\'evy spectral measure $M\in
\mathcal{M}(E)$, (6) gives that $M^{(\lambda)}\in \mathcal{M}(E)$.
Here is an extension for more general $\lambda$ but less general
$M$.

\begin{prop}
Suppose that $\lambda$ and $M$ are Borel measures on $(0,\infty)$
and $E$ respectively, such that
\begin{multline}
\int_{E\setminus\{0\}}\big[ ||x||\int_0^{||x||^{-1}}t\,
\lambda(dt) +
\lambda(\{s:s>||x||^{-1}\})]M(dx)= \\
\int_0^{\infty}\big[ t\int_{\{0< ||x||\le t^{-1}\}} ||x||M(dx)+
M(\{||x|>t^{-1}\})\big]\lambda(dt)<\infty.
\end{multline}
Then $M^{(\lambda)}$ and $M$ are L\'evy spectral measures on $E$
and $\lambda$ is a L\'evy spectral measure on $(0,\infty)$.
\end{prop}
\textit{Proof.} From Araujo-Gine (1980), Theorem 6.3, we know that
if a measure $M$ integrates $\min(1,||x||)$ on a Banach space $E$
then it is a L\'evy spectral measure. Condition (9) is just that
integrability condition for $M^{(\lambda)}$. [It was written in a
form of a sum of two integrals to indicate two different behaviors
of L\'evy measures: on open neighborhoods of zero and their
complements.]

Finally, from the inequality  $\min(1,t)\cdot\min(1,||x||)\le
\min(1, t||x||)$, for all $t \ge 0$ and $x\in E$, and formula (8)
we infer the remaining claims.

\begin{prop}
Suppose that $H$ is a real separable Hilbert space, and $M$  and
$\lambda$ are Borel measures on $H$ and on $(0,\infty)$,
respectively. Then $M^{(\lambda)}$ is a L\'evy spectral measure if
and only if
\begin{multline}
\int_{H\setminus\{0\}}\big[
||x||^2\int_0^{||x||^{-1}}t^2\,\lambda(dt) +
\lambda(\{s:s>||x||^{-1}\})]M(dx)= \\
\int_0^{\infty}\big[ t^2\int_{\{0< ||x||\le t^{-1}\}} ||x||^2
M(dx)+ M(\{||x|>t^{-1}\})\big]\lambda(dt)<\infty,
\end{multline}
Moreover, $M$ is a L\'evy spectral measures on $H$ and $\lambda$
is a L\'evy spectral measure on $\Rset$.
\end{prop}
\textit{Proof.} In a case of Hilbert space, for a measure $M$ to
be a L\'evy spectral measure it is necessary and sufficient that
$M$ integrates $\min(1,||x||^2)$; cf. Parthasarathy (1967),Chapter
VI, Theorem 4.10. The condition (10) is just the mentioned
integrability condition for $M^{(\lambda)}$. Furthermore, as
before we use the inequality $\min(1,t^2)\cdot\min(1,||x||^2)\le
\min(1, t^2||x||^2)$, for all $t \ge 0$ and $x\in H$, and the
proof is complete.

\medskip
\textbf{Examples. A).} Let $\emph{e}$ denotes the standard
exponential distribution with the density $e^{-t}, t>0$. Then
\begin{cor}
On any real separable Hilbert space H we have that
\begin{equation}
M^{(\emph{e})} \ \mbox{is a L\'evy spectral measure iff} \
\mbox{so is M}.
\end{equation}
\end{cor}
\textit{Proof.} Let $M$ be a L\'evy spectral measure on $H$. Since
we have that
\begin{multline*}
g(||x||):= ||x||^2\int_{0}^{||x||^{-1}}t^2 e^{-t}dt +
\int_{||x||^{-1}}^{\infty}e^{-t}dt =\\
2||x||^2 \big[1- e ^{-||x||^{-1}}(1+||x||^{-1})\big],
\end{multline*}
therefore $ g(||x||)\le 2 ||x||^2$, for $||x||\le1$. On the other
hand, using the power series expansion we also get $ \lim_{||x||
\to \infty} g(||x||)= 1$, which implies that $g(||x||) \le K$, for
$||x||> 1$. Consequently,
\begin{multline*}
\int_{H\setminus\{0\}}\min(1,||x||^2)M^{(\emph{e})}(dx) =
\int_{\{0<||x||\le1\}}g(||x||)M(dx)\\
+\int_{\{||x||>1\}}g(||x||)M(dx)\le \\
2\int_{\{0<||x||\le1\}}||x||^2M(dx) +K\int_{\{||x||>1\}}M(dx)<
\infty
\end{multline*}
and Proposition 2 gives that $M^{(\emph{e})}$ is spectral measure.

Converse implication also follows from Proposition 2, and thus the
proof is complete.

\medskip
\textbf{B).} Let $\rho_{\alpha}(dt):= t^{-\alpha-1}e^{-t}dt$ be a
measure on $(0,\infty)$.
\begin{cor}
On any Hilbert space $H$, if $M^{(\rho_{\alpha})}$ is a L\'evy
spectral measure then so is M and
\begin{multline*}
0< \alpha <2, \int_{||x||>1}||x||^{\alpha}M(dx)< \infty, \\
\mbox{ and for all $s>0$}, \ \ \int_{\{0<||x||\le
1\}}||x||^{\alpha}e^{-s/||x||}M(dx)< \infty.
\end{multline*}
Conversely, if $0<\alpha<2$ and
$\int_{\{||x||>0\}}||x||^{\alpha}M(dx)< \infty$ then
$M^{(\rho_{\alpha})}$ and $M$ are L\'evy spectral measures.
\end{cor}
\textit{Proof.} Let us introduce a function
\begin{multline*}
h_{\alpha}(||x||):= ||x||^2\int_{0}^{||x||^{-1}}t^{(2-\alpha)
-1}e^{-t}dt + \int_{||x||^{-1}}^{\infty}t^{- \alpha-1}e^{-t}dt \\
=||x||^{\alpha}\big[ \int_0^1t^{(2-\alpha) -1}e^{-t/||x||}dt
+\int_1^{\infty}t^{-\alpha -1}e^{-t/||x||}dt\big].
\end{multline*}
From Proposition 2, we have that $\int_{H\setminus
\{0\}}h_{\alpha}(||x||)M(dx)<\infty$. Consequently,
$h_{\alpha}(||x||) <\infty$ for $M-a.a.x$. Consequently, $2-\alpha
>0$ (from the first integral) and $\alpha>0$. Since we also have
that
\begin{multline}
\int_{H\setminus \{0\}}h_{\alpha}(||x||)M(dx)=
\int_0^1t^{(2-\alpha)
-1}\big(\int_{H\setminus\{0\}}||x||^{\alpha}e^{-t/||x||}M(dx)\big)dt \\
+\int_1^{\infty}t^{-\alpha
-1}\big(\int_{H\setminus\{0\}}||x||^{\alpha}e^{-t/||x||}M(dx)\big)dt
<\infty.
\end{multline}
Hence, the function $t \to
\int_{H\setminus\{0\}}||x||^{\alpha}e^{-t/||x||}M(dx) <\infty$ for
almost all (Lebesgue measure) $t \in \Rset^+$. Its monotonicity
gives that it is finite for all $t>0$. Finally note that
\[
e^{-1}\int_{\{||x||>1\}}||x||^{\alpha}M(dx)\le
\int_{\{||x||>1\}}||x||^{\alpha}e^{-1/||x||}M(dx) <\infty,
\]
which completes the proof.

The converse part follows from the fact that $h_{\alpha}(||x||)\le
2((2-\alpha)\alpha )^{-1}||x||^{\alpha}$ and Proposition 2. Thus
the proof is complete.

\medskip
\textbf{3. Random integral representations.} Random integral
representation method allows to represent a random variable, more
precisely its probability distribution, as a probability
distribution of random integrals of the form
$\int_{(a,b]}h(t)dY(r(t))$, where $Y$ is the L\'evy process
(process with stationary independent increments, cadlag paths and
starting from the origin), $h$ is a real deterministic and $r$ is
deterministic and monotone with positive values (deterministic
time change). A such method of a description of measures was
introduced in Jurek-Vervaat (1983) for selfdecomposable measures;
cf. also Jurek-Mason(1993), Chapter 3, Jurek (1982, 1985, 1988).

\medskip
\begin{thm}
Let $\lambda(\cdot)$ be a Borel measure on positive half-line that
is finite on sets bounded away from zero and let
$\Lambda(t):=\lambda(\{s>0: s>t\}),\, t>0$. Further, let
$Y(t),t\ge0$, be a cadlag L\'evy process with values in a Hilbert
space $H$ and $Y(1)$, as infinitely divisible random element, is
described by a triple $[a, R, M]$. Then in order that the limit
\begin{equation}
I_{(\alpha,\beta]}:=\int_{(\alpha,\,\beta]}\,t\,dY(\Lambda(t))\to
\int_{(0,\,\infty)}\,t\,dY(\Lambda(t))=: I_{(0,\infty)},
\end{equation}
exits in distribution ,as $\alpha \downarrow 0$ and $\beta
\uparrow \infty$, it is sufficient and necessary that
\begin{multline}
\int_0^{\infty}t \lambda(dt)< \infty,\ \mbox{provided} \ \ a \neq
0; \int_0^{\infty}t^2\lambda(dt)<\infty, \ \mbox{provided} \ \
R\neq 0;
\\
\int_{\{0<||x||\le 1\}}||x||\int_{||x||^{-1}}^{\infty} t
\lambda(dt)M(dx)+ \int_{\{||x||>1\}}||x||\int_0^{||x||^{-1}}t
\lambda(dt)M(dx)<\infty; \\
 \mbox{and} \ M^{(\lambda)} \ \mbox{is a L\'evy spectral measure.}
\end{multline}
Furthermore, if the limit $I_{(0,\infty)}$ has representation
$[a^{(\lambda)},R^{(\lambda)},M^{(\lambda)}]$ then
\begin{multline}
a^{(\lambda)}=(\int_0^{\infty}t \lambda(dt))\cdot a +
\int_0^{\infty}\int_{H\setminus\{0\}}[1_B(tx)-1_B(x)]\,t\,x\,
M(dx)\,\lambda(dt); \\
R^{(\lambda)}=(\int_0^{\infty}t^2 \lambda(dt))\cdot R;\ \quad \ \
M^{(\lambda)}(A)= \int_0^{\infty}\int_{H\setminus \{0\}}1_A
(tx)M(dx)\lambda(dt).
\end{multline}
\end{thm}
\textit{Proof}. From the definition of random integrals
$W_{(\alpha,\,\beta]}:=\int_{(\alpha,\beta]}h(t)dY(\tau(t))$,
where $h:(\alpha,\beta]\to \Rset$ and $\tau :(\alpha,\beta]\to
\infty$ are deterministic functions,
 and $\tau$ is monotone one, and $Y$ is a cadlag L\'evy process, we
have that
\begin{equation}
\mathbb{E}[e^{i<y,\,W_{(\alpha,\beta]}>}]=
\exp\int_{(\alpha,\beta]}\Big(\log
\mathbb{E}[e^{i<h(t)\,y,Y(1)>}]\Big) (\pm)\,d\tau(t),
\end{equation}
where one takes the  sign "+" for nondecreasing $\tau$ and "-" for
nonincreasing $\tau$; cf. Jurek-Vervaat(1983) or Jurek-Mason
(1993), Section 3.6. Hence using the L\'evy-Khintchine formula (2)
and taking $h(t)=t$ and $\tau=\Lambda$ in (16), we conclude that
random element $I_{(\alpha,\beta]}$ has an infinitely divisible
distribution with the triple $[a^{(\lambda)}_{(\alpha,\beta]},
R^{(\lambda)}_{(\alpha,\beta]}, M^{(\lambda)}_{(\alpha,\beta]}]$
given as follows
\begin{multline}
a^{(\lambda)}_{(\alpha,\beta]}=(\int_{(\alpha,\beta]}t
\lambda(dt))\cdot a + \int_{(\alpha,\beta]}t
\int_{H\setminus\{0\}}[1_{B}(tx)-1_{B}(x)] x
M(dx)\lambda(dt); \\
R^{(\lambda)}_{(\alpha,\beta]}=(\int_{(\alpha,\beta]}t^2 \lambda(dt))\cdot R;\\
M^{(\lambda)}_{(\alpha,\beta]}(A)=\int_{(\alpha,\beta]}\int_{H\setminus
\{0\}} 1_A(tx)\,M(dx)\,\lambda(dt)=M^{(\lambda
|_{(\alpha,\beta]})} (A),
\end{multline}
where the triple $[a,R,M]$ comes from the L\'evy-Khintchine
representation of the  infinitely divisible random element $Y(1)$.
As $\alpha \downarrow 0$ and $\beta \uparrow \infty$ then
$M^{(\lambda)}_{(\alpha,\beta]} \uparrow M^{(\lambda)}\in
\mathcal{M}(H)$, Gaussian covariance operators
$R^{(\lambda)}_{(\alpha,\beta]} \to R^{(\lambda)}$. Finally, for
the shift part note that
\[
|[1_B(tx)-1_B(x)]|=1 \ \mbox{iff}\ \ 1<||x||\le t^{-1} \ \
\mbox{or} \ \ t^{-1}<||x||\le 1,
\]
and the second  summand for a shift vector in (17) exits as a
Bochner integral on the product space $(0,\infty) \times
(H\setminus\{0\})$. Consequently, $I_{(\alpha,\beta]}$ converge in
distribution to $I_{(0,\infty)}$ by Parthasarathy (1968),Theorem
5.5,  because $M^{(\lambda)}_{(\alpha,\beta]} \uparrow
M^{(\lambda)}\in \mathcal{M}(H)$. Thus the proof is complete.

\begin{rem}
Since the random functions $\beta \to I_{(\alpha,\beta]}$ and
$\alpha \to I_{(\alpha,\beta]}$ have independent increments
(because so has L\'evy process $Y$) we infer that all three modes
of convergence: almost surly, in probability and in distribution)
are equivalent; cf. Araujo-Gine (1980), Chapter 3, Theorem 2.10,
p. 105.
\end{rem}

\medskip
As in previous papers Jurek\&Vervaat (1983) or Jurek (1982, 1985,
1988) here we introduce the following \emph{random integral
mapping}
\begin{equation}
\mathcal{K}^{(\lambda)}(\mu):=\mathcal{L}(\int_0^{\infty}t\,
dY_{\mu}(\Lambda(t)))\in ID
\end{equation}
where $Y_{\mu}(t), t \ge 0$ is a cadlag L\'evy process such that
$\mathcal{L}(Y_{\mu}(1))a=\mu$ and $\Lambda (\cdot)$ is the
cumulative distribution function or its tail function---note that
from Proposition 2, $\lambda$, as a L\'evy spectral measure, is
finite on any half-line $(a, \infty),\, a>0$.
\begin{cor}
For probability measures of the form
$\mathcal{K}^{(\lambda)}(\mu)$ one has
\[
\mathbf{E}[e^{i<y,\,\int_0^{\infty}t\,dY_{\mu}(\Lambda(t))>}]
 =
\exp \int_0^{\infty}\log \mathbf{E}[e^{it<y,Y_{\mu}(1)>}]
\lambda(dt),
\]
where $y\in E'$. Furthermore, the random integral mapping
$\mathcal{K}^{(\lambda)}(\mu)$ has the following algebraic
properties
\[
\mathcal{K}^{(\lambda)}(\mu_1 \star\mu_2)=
\mathcal{K}^{(\lambda)}(\mu_1)\star\mathcal{K}^{(\lambda)}(\mu_2),
\ \  \mathcal{K}^{(\lambda_1 + \lambda_2)}(\mu)=
\mathcal{K}^{(\lambda_1)}(\mu) \star
\mathcal{K}^{(\lambda_2)}(\mu)
\]
\end{cor}
\textit{Proof.} The first is consequence of the definition of
random integrals; cf. for analogous results in Jurek and Vervaat
(1983) or Jurek and Mason (1993), Lemma 3.6.4. The second equality
is a consequence of the formula (16) when one takes integrals over
positive half-line and the function $h(t)=t$.

One of the advantages of random integral representation is that it
allows easily to incorporate space and time changes. Here is an
example.
\begin{cor}
For $a\in\Rset,\,c>0$ and a random integral
$\int_{(\alpha,\beta]}h(t)dY(\tau(t))$, where $h:(\alpha,\beta]\to
\Rset$ and $\tau :(\alpha,\beta]\to \infty$ are deterministic
functions, and $\tau$ is monotone one, and $Y$ is a cadlag L\'evy
process, we have
\begin{multline}
\Big(\mathcal{L}\Big(
a\,\int_{(\alpha,\beta]}h(t)dY(\tau(t))\Big)^{\ast
c}\Big)^{\widehat{}}(y)=
\mathcal{L}\Big(\int_{(\alpha,\beta]}a\,h(t)dY(c\,\tau(t))\Big)(y)= \\
\exp\int_{(\alpha,\beta]}\log
\mathbf{E}[e^{i\,a\,h(t)<y,Y_{\mu}(c)>}]d\tau(t), \ \ y\in H.
\end{multline}
It is also true for integrals over half line, provided they exist.
\end{cor}
\textit{Proof.} Use (16) and the fact that $\mathcal{L}(Y(c))=
(\mathcal{L}(Y(1)))^{\ast c}$.

\medskip
\textbf{4. Two applications of the random integral method.}

\textbf{A)}. \emph{Free infinite divisibility.} As in Section 2,
let $\emph{e}(dt)$ denotes the standard exponential distribution.
Then from Example \textbf{A} we infer that $M^{(\emph{e})}$ is
L\'evy spectral measure (on $H$) whenever so is $M$. Furthermore
by Theorem 1, formula (14), $R^{(\emph{e})}= 2\,R$, and
\begin{multline}
a^{(\emph{e})}= a
+\int_{\{||x||>1\}}\,x(1-e^{-||x||^{-1}}(1+||x||^{-1}))M(dx) \\
+
\int_{\{0<||x||\le1\}}\,x\,e^{-||x||^{-1}}(1+||x||^{-1}))M(dx)
\end{multline}
exits in a Bochner sense. To this end note that $ \lim _{s\to
0}\,s\,e^{-s^{-1}}(1+s^{-1})= 0$  and $\lim _{s\to
\infty}\,s(1-e^{-s^{-1}}(1+s^{-1}))=0$.

Furthermore,  for $r>0$ and a Borel $D$ of the unit sphere
$S=\{x:||x||=1\}$, let us define \emph{L\'evy spectral function}
$L_M(D;r)$ associated with the measure $M$ as follows
\begin{equation}
L_M(D;r)= M(\{x: x\,||x||^{-1}\in D \ \mbox{and} \ \ ||x||>r \}).
\end{equation}
Then using (7) we  get
\begin{equation}
L_{M^{(e)}}(D;r)= \int_0^{\infty}\,L_M(D; rt^{-1})\,e^{-t}dt= r
\int_0^{\infty}\,L_M(D; s^{-1})\,e^{-r s}ds,\ \ r>0.
\end{equation}
Hence, $r^{-1}L_{M^{(e)}}(D;r), r>0,$ is a Laplace transform of
(unique) function $L_M(D;s^{-1})$ and thus $M^{(e)}$ uniquely
determines $M$. Hence, $a^{(e)}$ and $M$ uniquely identifies $a$.
All in all with Theorem 1 we conclude that
\begin{multline}
\mathcal{K}^{\emph{(e)}}: ID \ni \mu \to
\mathcal{L}(\int_0^{\infty}t\,dY_{\mu}(1-e^{-t}))\in \mathcal{E}:=
\mathcal{K}^{(e)}(ID) \\ \mbox {is well defined one-to-one random
integral mapping}, \ \ \ \  \
\end{multline}
where $Y_{\mu}(t), t\ge 0,$ is a cadlag L\'evy process such that
$\mathcal{L}(Y_{\mu}(1))=\mu$. Consequently, we obtained
convolution a subsemigroup $\mathcal{E} \subset ID$, which is
characterized among infinitely divisible by the triples
$[a^{(e)},R^{(e)},M^{(e)}]$ given by (15) and the kernel (2).

For a  probability measure $\nu$, let $\hat{\nu}$ denotes its
Fourier transform (characteristic function). In terms of Fourier
transforms elements representable as
$\mathcal{K}^{\emph{(e)}}(\cdot)$ are described as follows
\begin{cor}
In order for a function $g:H \to \Cset$ to be  a characteristic
function of a measure  from the convolution semigroup
$\mathcal{E}$ it is necessary and sufficient that
\begin{multline}
g(y)= \exp \big [i<y,a> - <y,Ry> +\\
\int_{H\setminus\{0\}}\Big(\frac{1}{1-i<y,x>} -1- i<y,x>
1_{\{||x||\le1\}}\Big)M(dx)],
\end{multline}
where $a \in H$, $R$ is non-negative, self-adjoint, trace operator
and $M$ is a Borel measure that integrates the function
$\min(1,||x||^2)$ over $H$. In fact, $g$ is a characteristic of
the measure $\mathcal{K}^{(e)}([a,R,M])$.
\end{cor}
One gets (24) by putting into (2) the triplet: the vector
$a^{(e)}$, the covariance  operator $R^{(e)}$ and the L\'evy
spectral measure $M^{(e)}$.
\begin{rem} One has two possibilities of looking at the class
$\mathcal{E}$. Either, as a subset of $ID$ with the triples
$[a^{(e)},R^{(e)},M^{(e)}]$ and the kernel $\Phi$ from formula (2)
or as a set of probability distributions given by triples [a,R,M]
and but with a new kernel
\begin{multline}
\Phi_1 (y):=\big [i<y,a> - <y,Ry> +\\
\int_{H\setminus\{0\}}\Big(\frac{1}{1-i<y,x>} -1- i<y,x>
1_{\{||x||\le1\}}\Big)M(dx)], \ y\in H.
\end{multline}
Note that both kernels are additive in $a, R$ and $M$, i.e., sums
of those parameters  correspond to convolution of probability
measures.

(Compare similar comments in  Jurek-Vervaat (1983) formula (4.3),
pages 254-255, for the L\'evy class $L$ of selfdecomposable
distributions.)
\end{rem}

\medskip
\begin{prop}
Let $I^{(e)}:= \int_0^{\infty}tdY(1-e^{-t})$ and
$\phi_{I^{(e)}}(y)$, and $\phi_{Y(1)}(y),\,\,y\in H$ are
characteristic functions of $I^{(e)}$ and $Y(1)$, respectively.
Then
\begin{multline}
\log
\phi_{I^{(e)}}(y)=\int_0^{\infty}\log\phi_{Y(1)}(ty)e^{-t}dt,\\
\log \phi_{Y(1)}(y)= \mathfrak{L}^{-1}[s^{-1}\log
\phi_{I^{(e)}}(s^{-1}y;x)]|_{x=1}, \ \ \ \
\end{multline}
where  for each $y\in H$, $\mathfrak{L}^{-1}$ is the inverse of
the Laplace transform of the function $s^{-1}\log
\phi_{I^{(e)}}(s^{-1}y)$. Hence, the mapping $\mathcal{K}^{(e)}:
ID(H)\to \mathcal{E}$ is an algebraic isomorphism between
convolution semigroups and for its inverse
$(\mathcal{K}^{(e)})^{-1}$ we have
\[
((\mathcal{K}^{(e)})^{-1}(\rho))^{\widehat{}}(y) = \exp
\mathfrak{L}^{-1}[s^{-1}\log \hat{\rho}(s^{-1}y;\,x)]|_{x=1},\ \ y
\in H.
\]
\end{prop}
\emph{Proof.} From Corollary 3 we have
\[
\log \phi_{I{(e)}}(y) = \int_0^{\infty}\log
\phi_{Y(1)}(ty)e^{-t}dt.
\]
Putting, for each fixed $y \in H$,
\[
f_y(s):=\log \phi_{I^{(e)}}(sy) \ \mbox{and} \ \ g_y(s):= \log
\phi_{Y(1)}(sy), \ \mbox{for} \ \ s\in \Rset,
\]
and using the above relation we get
\[
f_y(s)=\int_0^{\infty}g_y(ts)e^{-t}dt,\  \ f_y(s^{-1})= s
\int_0^{\infty}g_y(x)e^{-sx}dx, \ s>0.
\]
Consequently,
$\frac{1}{s}\,f_y(\frac{1}{s})=\mathfrak{L}[(g_y(x);s]$ is the
Laplace transform evaluated at $s$. This completes the proof.

\medskip
In an abstract semigroup $(\mathcal{G},\circ)$ and element $g\in
\mathcal{G}$ is said to be \emph{infinitely divisible} if for each
natural $n \ge 2$ there exists $g_n\in\mathcal{G}$ such that
n-times operation $g_n \circ g_n \circ ... \circ g_n = g$; cf.
Hilgert, Hoffman, Lawson (1989). By $ID(G)=(ID(G),\circ)$ we
denote a set of all $\circ$-infinite divisibility elements.

D. Voiculesu and others developed a theory of "a free
probability". For our needs here let us  recall briefly that with
a probability measure $\mu$, on a real line, one associates a
complex valued function $\emph{V}_{\mu}(z):=F^{-1}_{\mu}(z)-z,
z\in \mathcal{D}$, where $\mathcal{D}$ is an appropriately
selected domain in upper complex half-plane and
$F_{\mu}(z):=1/G_{\mu}(z)$, where
\begin{equation}
G_{\mu}(z):=\int_{-\infty}^{\infty}\frac{1}{z-t}\,\mu(dt)\, \ \
\mbox{is called the Cuchy transform}.
\end{equation}
For two probability measures, on the real line, $\mu$ and $\nu$,
if the sum $\emph{V}_{\mu}(z)+\emph{V}_{\nu}(z)$ corresponds to
another probability measure then we denote it by $\mu \Box \nu$.
Hence one can introduce $\Box$-infinite divisibility and a
semigroup $(ID(\mathcal{P}(\Rset),\Box)$. From Bercovici and Pata
(1999) and Barndorff-Nielsen and Thornbjornsen (2002) we have that
\begin{multline}
\nu \in ID(\mathcal{P}(\Rset),\Box)\ \ \mbox{iff}  \ \
z\,V_{\nu}(z^{-1})= az +\sigma^2z^2 +\\
\int_{\Rset\setminus\{0\}}\Big(\frac{1}{1 - zx} -1- zx
1_{\{|x|\le1\}}\Big)M(dx)], \ z\in\Cset^-, \ \ \ \
\end{multline}
where $a,\sigma \in \Rset$, and $M$ integrates $\min(1,|x|^2)$
over $\Rset$, i.e., $M$ is a L\'evy spectral measure on real line.
\begin{cor}
A probability measure $\nu$, on $\Rset$, is $\Box$-infinitely
divisible if and only if there exist a unique reals $a$ and
$\sigma^2$ and a L\'evy spectral measure $M$  such that
\begin{equation}
(it)\,V_{\nu}((it)^{-1})=\log(\mathcal{K}^{(e)}(\mu))^{\widehat{}}(t)=
\log \Big(\mathcal{L}
(\int_0^{\infty}sdY_{\mu}(1-e^{-s}))\Big)^{\widehat{}}(t), \ \
t\in\Rset,
\end{equation}
where $(Y_{\mu}(t), t\ge 0)$ is a L\'evy process such that
$\mathcal{L}(Y_{\mu}(1))=\mu =[a,\sigma^2,M]$. In other words,
functions $t\to e^{it\,V_{\nu}((it)^{-1})}, t\in\Rset,$ are
characteristic functions and a class of measures corresponding to
them coincides with the class $\mathcal{E}$. Furthermore, for the
Voiculescu transform we have $V_{\nu}(it)=
it\log(\mathcal{K}^{(e)}(\mu))^{\widehat{}}(-t^{-1}), \ \
t\in\Rset$.
\end{cor}
\emph{Proof.} Use Corollary 3 for $E=\Rset$ and then apply Theorem
1 with the formula (18).
\begin{rem}
Since both the kernel $\Phi_1$ given by (25), (that appeared in
the description of the class $\mathcal{E}$), and the kernel $\Phi$
given by (2), (that is the classical kernel from the
L\'evy-Khintchine), contain identical parameters $a,\,R$ and $M$
one has a natural identification between those convolution
semigroups. (See Remark 2). Sometimes it is called Bercovici-Pata
bijection between free and classical infinite divisibility. In the
approach presented here we have explicitly constructed the
isomorphism $\mathcal{K}^{(e)}$ between the semigroups
$\mathcal{E}$ and $ID$.
\end{rem}
\medskip
For the formula (24) in Corollary 5, or more precisely for the
existence of L\'evy process $Y$, it is necessary that $\mu$ is
$\ast$-infinitely divisible, while the Cauchy transform $G_{\mu}$,
given by (27), is defined for any (finite) measure. Here is a way
of avoiding that difficulty. For any finite measure $m$, on a
Hilbert or Banach space, let $e(m)$ denotes the compound Poisson
measure. Since it is $\ast$-infinitely divisible, (i.e., $e(m)\in
ID(H)$), we can insert it into a L\'evy (compound Poisson) process
$Y_{e(m)}(t), t\ge 0$. Consequently,
\begin{multline}
\mbox{if} \ \ \textbf{G}_m(y):=\int_{H}\frac{1}{1-i<y,x>}\,m(dx),
\ y \in H,\ \ \mbox{then}\\ \textbf{F}_m(y):=
\log(\mathcal{K}^{(e)}(e(m))^{\widehat{}}(y)=
\textbf{G}_m(y)-m(H)= \int_H\,\frac{i<y,x>}{1-i<y,x>}m(dx).
\end{multline}
To see that equalities recall that
$e(m)^{\widehat{}}(y)=\exp(\hat{m}(y)-m(H))$ and this with (18)
and Corollary 3 give the above formula.

\medskip
\textbf{B).} \emph{Tempered stable probability measures.} Let us
consider the example \textbf{B} from Section 2 for measures
$\rho_{\alpha}(dt)=t^{-\alpha-1}e^{-t}dt$ on positive half-line
with $0< \alpha <2$. Let us assume that
\begin{equation*}
\int_H \, ||x||^{\alpha}M(dx)<\infty \ \ \mbox{and $M$ is a Borel
 measure on $H$.}
\end{equation*}
In the sequel, by $ID_{\alpha}$ we denote those infinitely
divisible whose L\'evy spectral measures satisfy the above
integrability condition. Consequently, from Corollary 2 we have
that both $M^{(\rho_{\alpha})}$ and $M$ are L\'evy spectral
measures and from Theorem 1 we get $R^{(\rho_{\alpha})}=
\Gamma(2-\alpha)R$ is covariance operator of Gaussian measure. For
the shift vector $a^{(\rho_{\alpha})}$, we need three integrals;
cf. formula (14). Firstly, note that
\begin{multline*}
\int_{\{0<||x||\le 1\}}||x|| \int_{||x||^{-1}}^{\infty}
t^{-\alpha} e^{-t}dt\,M(dx) \le \\ \int_{\{0<||x||\le
1\}}||x||^2M(dx)\int_{1}^{\infty} t^{(2-\alpha)-1}
e^{-t}dt\,M(dx)< \infty, \ \mbox{for} \ \ 0< \alpha <2.
\end{multline*}
And secondly, note that
\begin{multline*}
 \int_{\{||x||> 1\}}||x||\int_0^{||x||^{-1}}t^{-\alpha} e^{-t}dtM(dx)=
\int_{\{||x||> 1\}} ||x||\gamma(1-\alpha, ||x||^{-1})M(dx)\\
= \int_{\{||x||>1\}}||x||^{\alpha}\Big(\sum_{n=0}^{\infty}
\frac{(-1)^n}{n!(\alpha+n)}
\cdot\frac{1}{||x||^{n}}\Big)M(dx)<\infty,\ \mbox{for} \
0<\alpha<1.
\end{multline*}
Consequently, for $0<\alpha<1$, the random integral mapping
\begin{multline}
\mathcal{K}^{(\rho_{\alpha})}: ID_{\alpha} \ni \mu \to
\mathcal{L}(\int_0^{\infty}t\,dY_{\mu}(\Gamma(\alpha,t))) \in
\mathcal{TS}_{\alpha}:= \mathcal{K}^{(\rho_{\alpha})}(ID_{\alpha})
\end{multline}
is well defined. In above $\gamma(\alpha,x)$ and $\Gamma(\alpha,
x)$ denote the incomplete Euler's gamma functions, i.e.,
\[
\gamma(\alpha,x)=\int_0^x\,t^{\alpha-1}e^{-t}dt, \ x>0,\,
(\Re\alpha>0); \ \ \Gamma(\alpha,
x)=\int_x^{\infty},t^{\alpha-1}e^{-t}dt, \ x>0.
\]
 Following Rosi\'nski(2002) measures from the class
$\mathcal{TS}_{\alpha}$ are called \emph{tempered stable
distributions}. In fact, they were introduced as infinitely
divisible measures (on $\Rset^d$) without Gaussian parts with
spectral measures of the form $M^{(\rho_{\alpha})}$ from Example
\textbf{B} in Section 2. [In Rosi\'nski (2004),
$M^{(\rho_{\alpha})}$ appears as Lemma 2.2.] Let us mention here
that tempered stable processes are of importance in statistical
physics as they exhibits different local and global behavior; (cf.
Mantegna and Stanley (1994), Novikov (1995), Kaponen (1995); comp.
Corollary 7 below. However, our point of interest is that the
tempered stable probability measures admit random integral
representation as well.
\begin{prop}
Assume that $0<\alpha<1$. Let $I^{(\rho_{\alpha})}:=
\int_0^{\infty}tdY(\Gamma(-\alpha,t))$ and
$\phi_{I^{(\rho_{\alpha}}}(y)$, and $\phi_{Y(1)}(y),\,\,y\in H$,
are characteristic functions of $I^{(\rho_{\alpha})}$ and $Y(1)$,
respectively. Then
\begin{multline}
\log
\phi_{I^{(\rho_{\alpha})}}(y)=\int_0^{\infty}\log\phi_{Y(1)}(ty)t^{-\alpha-1}e^{-t}dt,\\
\log \phi_{Y(1)}(y)= \mathfrak{L}^{-1}[s^{\alpha}\log
\phi_{I^{(\rho_{\alpha})}}(s^{-1}y;x)]|_{x=1}, \ \ \ \
\end{multline}
where  for each $y\in H$, $\mathfrak{L}^{-1}$ is the inverse of
the Laplace transform

$s^{\alpha}\log \phi_{I^{(\rho_{\alpha})}}(s^{-1}y)$. Hence, the
mapping $\mathcal{K}^{(\rho_{\alpha)}}: ID(H_{\alpha}) \to
\mathcal{TS}_{\alpha}$ is an algebraic isomorphism between
convolution semigroups. For its inverse
$(\mathcal{K}^{(\rho_{\alpha})})^{-1}$ we have
\[
((\mathcal{K}^{(\nu_{\alpha})})^{-1}(\rho))^{\widehat{}}(y) = \exp
\mathfrak{L}^{-1}[s^{\alpha}\log \hat{\nu}(s^{-1}y;\,x)]|_{x=1},\
\ y \in H, \ \ (\nu\in\mathcal{TS}_{\alpha}).
\]
\end{prop}
Proof is analogous to that of Proposition 3.
\medskip
\begin{rem} The previous result is also true for $1\le \alpha<2$
if one considers only L\'evy processes with symmetric spectral
measures $M$ and zero shifts $a$.
\end{rem}
Here is an example of usefulness of the random integral
representation of random variables. The result below is announced
in Rosi\'nski (2002); also cf. Rosi\'nski (2004). The proof  below
is based on the random integral representation of tempered stable
probability measures.
\begin{cor}
Let $0<\alpha<1$ and $X:=\int_0^{\infty}t\,dY(\Gamma(-\alpha,t)$
be $\Rset^d$-valued random vector with
$\mathcal{L}(Y(1))=[0,0,M]\in ID_{\alpha}$. Then
\[
(\mathcal{L}(s^{-1/\alpha}\,X)^{\ast s})\Rightarrow \eta_{\alpha},
\ \mbox{as}\ \ s\to 0,
\]
where $\eta_{\alpha}$ denotes the strictly stable law with
exponent $\alpha$
\end{cor}
\emph{Proof.} Using Corollary 4 with $a=s^{-1/\alpha}$ and $c=s$
we have
\begin{multline*}
((\mathcal{L}(\frac{1}{s^{1/\alpha}}X)^{\ast
s}))^{\widehat{}}(y)=\exp\int_0^{\infty}\,s\,\log \mathbf{E}[\exp
i\,t/s^{1/\alpha}<y,Y(1)>]t^{-\alpha-1}e^{-t}dt\\
=\exp\int_0^{\infty}\log\mathbf{E}[\exp i
u<y,Y(1)>]u^{-\alpha-1}e^{-u\,s^{1/\alpha}}du \ \ \mbox{(as $s\to 0)$} \\
\to \exp\int_{\Rset^d\setminus{\{0\}}}\int_0^{\infty}
[e^{iu<y,x>}-1-iu<y,x>1_{||x||\le1}(x)]u^{-\alpha-1}du\,M(dx)=\\
\exp[-c_{\alpha}\int_{\Rset^d\setminus{\{0\}}}|<y,x>|^{\alpha}(1-i\,\tan(\pi\alpha/2)\,
sign<y,x>)M(dx)],
\end{multline*}
where $c_{\alpha}>0$. (The last equality is obtained via contour
integration; see any book on stable laws.) Finally, the last
formula is the characteristic function of a strictly stable
probability measure on $\Rset^d$.

\medskip
\noindent Institute of Mathematics, University of Wroc\l aw,
50-384 Wroc\l aw, Poland. [E-mail: zjjurek@math.uni.wroc.pl]
\end{document}